\documentclass[11pt]{amsart}
\usepackage{amsfonts,amssymb,latexsym,amsmath, amsxtra}
\usepackage[all]{xy}
\usepackage[dvips]{graphics}

\usepackage[OT2,T1]{fontenc}
\DeclareSymbolFont{cyrletters}{OT2}{wncyr}{m}{n}
\DeclareMathSymbol{\Sha}{\mathalpha}{cyrletters}{"58}

\theoremstyle{plain}
\newtheorem{theorem}{Theorem}[section]
\newtheorem{corollary}[theorem]{Corollary}

\newtheorem{lemma}[theorem]{Lemma}

\newtheorem*{conjecture*}{Conjecture}
\newtheorem*{refconjecture}{Refined Conjecture}
\theoremstyle{definition}

\theoremstyle{remark}
\newtheorem*{remark}{Remark}

\numberwithin{equation}{section}

\newcommand{\Z}{\mathbb Z}

\newcommand{\Q}{{\mathbb Q}}

\def\calM{\mathcal{M}}
\def\calN{\mathcal{N}}

\def\({\left(}
\def\){\right)}

\newcommand{\rank}{\text{rank}}
\newcommand{\crank}{\text{crank}}

\def\k2{\frac{k}{2}}

\begin{document}

\title{Asymptotics for Rank and Crank Moments}
\author{Kathrin Bringmann}
\address{Mathematical Institute\\University of
Cologne\\ Weyertal 86-90 \\ 50931 Cologne \\Germany}
\email{kbringma@math.uni-koeln.de}
\author{Karl Mahlburg}
\address{Department of Mathematics \\
Massachusetts Institute of Technology \\
 MA 02139\\ U.S.A.}
\email{mahlburg@math.mit.edu}
\author{Robert C. Rhoades}
\address{Ecole Polytechnique Fed\'eral\'e de Lausanne, 1015 Lausanne, Switzerland}
\email{robert.rhoades@epfl.ch}

\thanks {The first author was partially supported by NSF grant DMS-0757907. The second author was partially supported by NSA Grant 6917958.}

\date{\today}
\thispagestyle{empty} \vspace{.5cm}

\begin{abstract}
Moments of the partition rank and crank statistics have been studied for their connections to combinatorial objects such as Durfee symbols, as well as for their connections to harmonic Maass forms.  This paper proves a conjecture due to Bringmann and Mahlburg that refined a conjecture of Garvan.  Garvan's conjecture states that the moments of the crank function are always larger than the moments of the rank function, even though the moments have the same main asymptotic term.  The proof uses the Hardy-Ramanujan method to provide precise asymptotic estimates for rank and crank moments and their differences.
\end{abstract}

\maketitle

\section{Introduction and statement of results}

The theory of partitions has long motivated the study of hypergeometric series and automorphic forms.
A foundational example for the interplay between these fields is Euler's partition function
$p(n)$, which has the generating function
\begin{equation} 
P(q):= \sum_{n=0}^{\infty} p(n) q^{n} = \prod_{n=1}^{\infty}
\frac{1}{1-q^{n}};\end{equation} this coincides with $q^{\frac{1}{24}}/\eta(z)$ where $\eta(z):= q^{1/24}\prod_{n=1}^{\infty}(1-q^{n})$ is Dedekind's weight $1/2$ modular form.    On the other hand,
many partition generating functions do not naturally appear in the theory of modular forms but rather in more general automorphic contexts, such as the theory of harmonic Maass forms.
Ramanujan's famous {\it mock theta functions} are prime examples of this phenomenon \cite{Wat}.
As a result, the study of generating functions from partition theory has
inspired a number of important results about  mock theta functions
and harmonic Maass forms \cite{B08, BGM, BO2, BR}.
Harmonic  Maass forms are real analytic generalizations of  modular forms, 
in that they   satisfy    the same transformation law and (weak) growth conditions at 
cusps, but instead of being holomorphic functions of the complex upper half plane, they are only required to be annihilated by the weight $k$ hyperbolic Laplacian.

Many of these modern results originated in Ramanujan's original results on the arithmetic of the partition function \cite{Ram19, Ram21}.  Most famously, the three ``Ramanujan congruences'' state that for all $n$,
\begin{align}\label{eqn:Rama}
p(5n+4) & \equiv 0 \pmod{5} \\
p(7n+5) & \equiv 0 \pmod{7}  \notag \\
p(11n+6) & \equiv 0 \pmod{11}. \notag
\end{align}

In an effort to provide a combinatorial explanation of Ramanujan's congruences, Dyson introduced \cite{Dy} the {\it rank} of a partition, which is defined as
\begin{equation*}
\rank(\lambda) := \text{largest part of} \; \lambda - \text{number of parts of} \; \lambda.
\end{equation*}
He conjectured that the partitions of $5n+4$ (resp. $7n+5$)
form 5 (resp. 7) groups of equal size when sorted by their rank modulo $5$ (resp. 7).  Building on Dyson's observations, Atkin and Swinnerton-Dyer later proved Dyson's rank conjectures \cite{AS}.

Dyson further conjectured the existence of an analogous statistic, the {\it crank}, that would explain all three congruences simultaneously.  Garvan finally found the crank while studying $q$-series of the sort seen in Ramanujan's ``Lost Notebook'' \cite{g2}, and together with Andrews presented the following definition \cite{AndG}.  Let $o(\lambda)$ denote the number of ones in a partition, and define $\mu(\lambda)$ as the number of parts strictly larger than $o(\lambda).$  Then
\begin{equation}
\crank(\lambda) := \begin{cases} \text{largest part of} \; \lambda \qquad & \text{if} \; o(\lambda) = 0, \\
\mu(\lambda) - o(\lambda) & \text{if} \; o(\lambda) > 0.
\end{cases}
\end{equation}
Works of the first and second authors show that both the rank and crank also play key roles in understanding the infinitely many other congruences for $p(n)$ \cite{B0, Mah}.

It is more useful here to work with generating functions than combinatorial definitions.  If $\calM(m,n)$ and $\calN(m,n)$ are the number of partitions of $n$ with crank and rank $m$, respectively, then, aside from the anomalous case of $\calM(m,n)$ when $n=1$, the two-parameter generating functions may be written as  \cite{AndG, AS}
\begin{align}
C(x;q) & := \sum_{\substack{m \in \Z \\ n \geq 0}} \calM(m,n) x^m q^n =
\prod_{n \geq 1} \frac{1-q^n}{(1-xq^n)(1-x^{-1}q^n)} =
\frac{1-x}{(q)_\infty} \sum_{n \in \Z} \frac{(-1)^{n} q^{n(n+1)/2}}{1-xq^n},
\notag \\
R(x;q) & := \sum_{\substack{m \in \Z \\ n \geq 0}} \calN(m,n) x^m q^n =
\sum_{n \geq 0} \frac{q^{n^2}}{(xq;q)_n (x^{-1}q;q)_n} =
\frac{1-x}{(q)_\infty} \sum_{n \in \Z} \frac{(-1)^{n} q^{n(3n+1)/2}}{1-xq^n}.
\notag
\end{align}
Although the final expressions for $C(x;q)$ and $R(x;q)$ appear quite similar, their analytic behaviors are markedly different.  For example, if $x \neq 1$ is a fixed root of unity, then $C(x;q)$ is essentially a meromorphic modular form \cite{Mah}, whereas $R(x;q)$ corresponds to the holomorphic part of a harmonic Maass form \cite{BO2}.

In addition to the importance of the rank and crank statistics in the study of Ramanujan's
congruences,
a number of recent works highlight the importance of the  weighted
moments  of  the crank and rank statistics.
For example, Andrews obtained an elegant description of the smallest parts partition function $spt(n)$, the number of smallest parts among the integer partitions of $n$, in terms of $p(n)$ and $N_{2}(n): = \sum_{m} m^{2} N(m,n)$.  Andrews proved \cite{And08} that $$spt(n) = np(n) - \frac{1}{2} N_{2}(n).$$
Additionally, Andrews introduced {\it Durfee symbols} and some refinements in \cite{And07}, which rely heavily on the intrinsic combinatorial connection between Durfee symbols and crank and rank moments.  In joint work with Garvan, the first and second authors proved a variety of new Ramanujan-type congruences for ``marked'' Durfee symbols \cite{BGM}.
Finally, Atkin and Garvan's original study of rank and crank moments also had important applications to partition congruences \cite{AG}.

For an integer $k$, define the $k$-th {\it crank} (resp. {\it rank}) {\it moment} as
\begin{align}
M_k(n)&:=\sum_{m \in \Z} m^k\, \calM(m,n), \\
N_k(n)&:=\sum_{m \in \Z} m^k\, \calN(m,n). \notag
\end{align}
Both the crank and rank moments vanish when $k$ is odd due to the symmetries of the statistics \cite{AG}.
Let
\begin{align*}
C_{k}(q) := &\sum_{n} M_{k}(n)q^{n}\\
R_{k}(q) :=& \sum_{n} N_{k}(n) q^{n}
\end{align*}
be the generating functions for $M_{k}(n)$ and $N_{k}(n)$. As the generating functions $C(x;q)$ and $R(x;q)$ exhibit very different analytic properties, so do these generating functions.  The first two authors along with Garvan \cite{B, BGM} showed that while $C_{k}(q)$ is
essentially a quasimodular form, $R_{k}(q)$ is a {\it quasimock theta function}, defined as the holomorphic part of sums of weak Maass forms and their derivatives.

In this paper we focus on the unpublished observations of Garvan regarding the relative size of the crank and rank moments.
\begin{conjecture*}[Garvan]
Let $k \geq 1$ be an integer.
\begin{enumerate}
\item
As $n \rightarrow \infty,$ we have $M_{2k}(n) \sim N_{2k}(n)$.
\item
For all $n \geq 2,$ we have $M_{2k}(n) > N_{2k}(n).$
\end{enumerate}
\end{conjecture*}
\begin{remark}
Garvan's conjecture can be interpreted as stating that the distribution of the crank statistic is slightly ``wider'' than that of the rank, but not enough to affect the main asymptotic behavior.  This is unexpected, as there is little about the combinatorial definitions of the crank and rank that suggest any close relations.
\end{remark}

A refined conjecture was given by the first two authors in \cite{BM}.
\begin{refconjecture}[Bringmann-Mahlburg]
\label{MainConj}  Suppose that $k \geq 1$.
\begin{enumerate}
\item
As $n \rightarrow \infty,$
\begin{equation*}
M_{2k}(n) \sim N_{2k}(n) \sim \alpha_{2k} \cdot n^{k} \: p(n),
\end{equation*}
where $\alpha_k \in \Q$ is non-negative. 
\item
Garvan's inequality holds for all $n$, and as $n \rightarrow \infty,$
\begin{equation*}
D_{2k}(n):= M_{2k}(n)-N_{2k}(n) \sim \beta_{2k} \, n^{k-\frac{1}{2}} \cdot p(n),
\end{equation*}
where $\beta_{2k} \in \frac{\sqrt{6}}{\pi} \: \Q$ is positive
\end{enumerate}
\end{refconjecture}

\begin{remark}
The second part of the conjecture is not a refinement of Garvan's conjecture,
except for sufficiently large $n$.
\end{remark}
The first two authors proved this conjecture in the cases $k = 2$ and $4$, and outlined a general procedure for calculating the main asymptotic terms of all crank and rank moments.  However, this procedure did not generically imply that the main terms were equal, as the recurrences that relate cranks and ranks are quite complicated.  

In this paper, we prove the Bringmann-Mahlburg conjecture
by computing explicit asymptotic expansions for the rank and
crank moments.  We adopt the standard notation $B_n(x)$ for the Bernoulli polynomials in the statement of the following result.  

\begin{theorem}\label{thm:crank}
We have
$$
M_{2k}(n) =
 \pi \xi_{2k} (24 n-1)^{k - 3/4} I_{3/2}(y_{n}) + \widetilde{\xi}_{2k}
(24n-1)^{k-5/4}
I_{1/2}(y_{n})
 +O\(n^{k-7/4}\cdot n^{-1/4} \exp( y_{n}) \),
 $$
where
\begin{equation*}
y_{n}:= \frac{\pi}{6} \sqrt{24n-1}.
\end{equation*}
The constants are given by
\begin{equation*}
 \xi_{2k}:= (-1)^{k} 2 B_{2k}(1/2) \qquad \text{and} \qquad
\widetilde{\xi}_{2k}:=
-3 (2k)(2k-3)  \xi_{2k}+ \xi_{2k}',
\end{equation*}
where $\xi_0':=0$, and for $k>0$,
$$
\xi'_{2k}:=-\frac14 (2k)(2k-1)\xi_{2k-2}.
$$
\end{theorem}
\begin{remark} At $x = 1/2$, the Bernoulli polynomials evaluate to 
$B_{2k}(1/2) = \left(2^{1-2k}-1 \right)B_{2k},$ where $B_{2k}$ is the usual Bernoulli number.  In particular, we have that $\xi_{2k} = 2(-1)^k B_{2k}(1/2) >0$.
\end{remark}
\begin{remark}
Regardless of index, each Bessel function has the main asymptotic term $\displaystyle I_a(y) \sim \frac{e^y}{y}$ \cite{AAR}.
\end{remark}

Theorem \ref{thm:crank} follows from a circle method argument and
a recursive relation for the crank moment
generating functions given in \cite{AG}.    The key ingredient in studying the rank
moment generating functions is Atkin-Garvan's ``Rank-Crank PDE'' \cite{AG}, which we state precisely in Section \ref{sec:Rank}.   This PDE is a recursive formula
for the rank moment generating functions that involves triple products of crank moment generating functions.
We use an asymptotic expansion for the coefficients of products of crank moment generating functions (see Section \ref{sec:CrankProducts}) in order to deduce the following
rank moment asymptotics.
\begin{theorem}\label{thm:rank} For $k\ge 0$ we have
\begin{align*}N_{2k}(n) =& \pi \lambda_{2k} (24 n-1)^{k - 3/4} I_{3/2}(y_{n}) + \widetilde{\lambda}_{2k} (24n-1)^{k-5/4}I_{1/2}(y_{n})+O\(n^{k-7/4}\cdot n^{-1/4} \exp( y_{n}) \),
\end{align*}
where $\lambda_{2k} := \xi_{2k}$ and $$\widetilde{\lambda}_{2k}
: = \begin{cases} -3\cdot 2k(2k-3) \xi_{2k}
 -\frac34 2k(2k-1) \xi_{2k-2}& k >0, \\
  0 & k =0. \end{cases}
$$
\end{theorem}

Combining Theorems \ref{thm:crank} and \ref{thm:rank}, we obtain an asymptotic expansion for
the difference of the rank and crank moments.
\begin{corollary}
 For any $k \ge 1$ as $n\to \infty$ we have
$$D_{2k}(n) \sim  \frac12\cdot 2k(2k-1) \xi_{2k-2}  (24n-1)^{k-5/4}I_{1/2}(y_{n}) . $$
\end{corollary}
\noindent Finally, the constants in the refined conjecture can now be computed using the fact that 
$$p(n) \sim 2 \pi (24n-1)^{-3/4} I_{3/2}(y_{n}).$$
\begin{corollary}
The Bringmann-Mahlburg conjecture is true, with constants
\begin{align*}
\alpha_{2k} & = (-24)^{k} B_{2k}(1/2), \\
\beta_{2k} & =  \frac{\sqrt{6}}{\pi} \cdot 2k(2k-1)(-24)^{k-1} B_{2k-2}(1/2).
\end{align*}
\end{corollary}

\section{Crank Asymptotics and the proof of Theorem \ref{thm:crank}}
Here we prove a  modified version of Theorem \ref{thm:crank}, which reflects the Bessel function indices that arise most naturally when using the circle method.
 \begin{theorem}\label{thm:crank2}
 For $k\ge 0$ we have
\begin{align*}
M_{2k}(n) =& \pi \xi_{2k}\cdot (24 n-1)^{k - 3/4} I_{3/2-2k}(y_{n}) + \xi'_{2k} \cdot(24n-1)^{k-5/4}I_{3/2-2k+1}(y_{n})\\ & +O\(n^{k-7/4}\cdot n^{-1/4} \exp( y_{n}) \).
\end{align*}
  \end{theorem}
\noindent Theorem \ref{thm:crank} follows from Theorem \ref{thm:crank2} through a simple formula for shifting the indices of Bessel functions.
  \begin{lemma}\label{lem:besselFunc}
  For $\ell \in \Z$, we have the relation
 $$
 I_{3/2 - 2\ell} (y_{n}) = I_{3/2}(y_{n}) - \frac{3}{\pi} (24n-1)^{-1/2} (2\ell)(2\ell -3) I_{1/2}(y_{n}) + 
 O\left(n^{-1} I_{-1/2}(y_{n}) \right).
 $$
 \end{lemma}
 \begin{proof}
 The total shift is the result of successive applications of the Bessel function relation \cite{AAR}
 $$
 I_{a-1}(x)
 = \frac{2a}{x}I_{a}(x)
 + I_{a+1}(x).
 $$
 \end{proof}
 \begin{proof}[Proof of Theorem \ref{thm:crank2}]
The idea of the proof is to find a recursive formula for the leading order constants, which we then solve explicitly with Bernoulli polynomials, thus obtaining formulas for the constants $\xi_{2k}$ and $\xi'_{2k}$.
We begin with Atkin and Garvan's recurrence for the
crank moment generating functions in terms of divisor sums, found as equation (4.6) in \cite{AG}, namely
\begin{equation}\label{eqn:crankREC}
C_{2k}(q) = 2 \sum_{j=1}^{k} \binom{2k-1}{2j-1} \Phi_{2j-1}(q)C_{2k-2j}(q).
\end{equation}
Here we have denoted the $j$-th divisor function by
\begin{equation*}
\Phi_{j}(q):= \sum_{n=1}^{\infty} \sigma_{j}(n) q^{n},
\end{equation*}
where $\sigma_{j}(n) := \sum_{d\mid n} d^{j}$ is the $j$-th divisor sum.  These functions can be written in terms of the classical
Eisenstein series $E_{k}(z) := 1- \frac{2k}{B_{k}} \Phi_{k-1}(q),$ as
\begin{equation*}
\Phi_{k-1}(q) = -\frac{B_k}{2k} \cdot \left(E_k(q) - 1\right).
\end{equation*}

We next use the Hardy-Littlewood circle method, amplifying the arguments used in \cite{BM}.
Let $q:=e^{ -  2 \pi z}$ with Re$(z)>0$, $q_1:=e^{-\frac{2 \pi }{z} }$.
Then we have the transformation laws
\begin{align}
\label{EkE2}
E_k(q) &= (iz)^{-k} E_k(q_1) \qquad  \text{if} \; k>2, \\
E_2(q)&= (iz)^{ -2} E_2(q_1) +\frac{6}{ \pi z}. \notag
\end{align}
In general, suppose that we wish to estimate the coefficients in an expression of the form
\begin{equation}
\label{Pqz}
\sum_n a(n) q^n = c \, P(q) g(q_1) z^{-k} + \dots,
\end{equation}
where $c$ is a constant and $g(q)$ has a holomorphic $q$-series expansion
$$
\displaystyle
g(q)=1 +\sum_{ n>0} b(n)q^{ n}.
$$
Then the asymptotic contribution to $a(n)$ due to the term displayed on the right side of (\ref{Pqz}) is
\begin{equation*}
\label{BesselAsymptotic2}
c\cdot 2 \pi (24n-1)^{ \frac{k}{2}-\frac34  } I_{\frac32 -k}(y_n).
\end{equation*}

This immediately implies the following recursive result.
\begin{lemma}\label{lem:Eisen}
Suppose that $G(q) = \sum c(n)q^n = P(q)\widetilde{E}_{2k}(q),$ where $\widetilde{E}_{2k}(q) := E_{2a_1}(q)\cdots E_{2a_r}(q)$ has total weight $2k = 2a_1 + \dots + 2a_r.$  Then the coefficients $c(n)$ have an asymptotic expansion of the form
\begin{equation*}
c(n) =\pi  \alpha (24n-1)^{k - 3/4}I_{3/2 - 2k}(y_{n}) + \alpha' (24n-1)^{k-5/4}I_{5/2-2t}(y_{n}) + O\left(n^{k - 2} e^{y_{n}}\right) 
\end{equation*}
for some constants $\alpha$ and $\alpha'$. 
Furthermore, the coefficients of $G(q) \cdot E_{2t}(q) =: \sum c_{2t}(n)q^{n}$ satisfy
\begin{align*}
c_{2t}(n) & =\pi  \alpha (-1)^{t} (24n-1)^{k-3/4 +t}I_{3/2 - 2k - 2t}(y_{n})  \\
& + (\alpha'(-1)^{t} + \delta_{t=1} \cdot 6
\alpha) (24n-1)^{k-5/4+t}I_{5/2 - 2k-2t}(y_{n}) + O\left(n^{k- 2 +t} e^{y_{n}}\right).
\end{align*}
\end{lemma}

We apply this lemma to the crank moments by first ``unwinding'' (\ref{eqn:crankREC}) to obtain the formula
\begin{equation}\label{eqn:crankUnwind}
C_{2k}(q) = 2P(q) \cdot \sum_{a_{1}+2a_{2} +  \cdots + ka_{k} = k} \alpha_{a_{1}, a_{2}, \cdots, a_{k}} \Phi_{1}(q)^{a_{1}}\Phi_{3}(q)^{a_{2}}\cdots\Phi_{2k-1}(q)^{a_{k}}
\end{equation}
for some integer constants $\alpha$.  After rewriting in terms of Eisenstein series, Lemma \ref{lem:Eisen} now applies to each term in (\ref{eqn:crankUnwind}), which implies the existence of the asymptotic expansion in Theorem \ref{thm:crank2}.
Furthermore, (\ref{eqn:crankREC}) and Lemma \ref{lem:Eisen} then also immediately imply the following recurrence:
\begin{equation}\label{eqn:xiREC}
\xi_{2k} = 2 \sum_{j=1}^{k} \binom{2k-1}{2j-1} \frac{B_{2j}}{4j}(-1)^{j+1} \xi_{2k-2j}.
\end{equation}
Similarly, we may deduce the following recurrence for $\xi_{2k}'$,
\begin{align}\label{eqn:xi'REC}
\xi_{2k}' =& 2\sum_{j=1}^{k-1} \binom{2k-1}{2j-1} \frac{B_{2j}}{4j} (-1)^{j+1} \xi_{2k-2j}' - 2(2k-1)\frac{B_{2}}{4} \, 6\,  \xi_{2k-2} \\
=& 2\sum_{j=1}^{k-1} \binom{2k-1}{2j-1} \frac{B_{2j}}{4j} (-1)^{j+1} \xi_{2k-2j}' - \frac{(2k-1)}{2} \xi_{2k-2}. \notag
\end{align}

We solve both of these recurrences by using the Bernoulli polynomial identity
\begin{equation}\label{eqn:bern}
\sum_{j=1}^{k} \binom{2k-1}{2j-1} \frac{B_{2j}}{4j} B_{2k-2j}(1/2) = -\frac{1}{2} B_{2k}(1/2).
\end{equation}
This is a specialization of the general convolution sum
\begin{equation*}
\sum_{j=0}^n \binom{n}{j} B_j(x) B_{n-j}(y) = -(n-1)B_n(x+y) + n(x+y-1)B_{n-1}(x+y),
\end{equation*}
which can be found (along with other relevant formulas) in \cite{Dil}.  The fact that $B_n(1/2) = 0$ for all odd $n$ implies that our formula (\ref{eqn:bern}) is equivalent to the case $x = 0, y = 1/2.$  Applying induction to (\ref{eqn:xiREC}) and (\ref{eqn:bern}) and using the base case $\xi_0 = 2$ gives the claimed formula for $\xi_{2k}$.

 To obtain the formula for $\xi_{2k}'$, note that $\binom{2k-1}{2j-1} = \frac{(k-1)(2k-1)}{(k-j)(2k-2j-1)} \binom{2k-3}{2j-1}$.  Now the recurrence (\ref{eqn:xiREC}) inductively implies that the claimed formula for $\xi_{2k}'$ is correct, as the right side of (\ref{eqn:xi'REC}) evaluates to
\begin{align*}
&
 -\frac12  (2k-1)(k-1) \sum_{j=1}^{k-1} 2 
 \binom{2k-3}{2j-1} \frac{B_{2j}}{4j}  (-1)^{j+1} \xi_{2(k-1)-2j}
 - \frac{(2k-1)}{2} \xi_{2k-2}    \\
=&     -\frac12  (2k-1)(k-1) \xi_{2(k-1)} - \frac{(2k-1)}{2} \xi_{2k-2}
= -\frac12k(2k-1)\xi_{2k-2}
= \xi_{2k}'.
\end{align*}
\end{proof}
We note for later purposes the identity
\begin{equation}
\label{eqn:xi'ID}
 2\sum_{j=1}^{k-1} \binom{2k-1}{2j-1} \frac{B_{2j}}{4j} (-1)^{j+1} \xi_{2k-2j}'
 = \left(1-\frac{1}{k} \right)\xi_{2k}'.
\end{equation}

\section{Asymptotics for products of crank moments}\label{sec:CrankProducts}
The rank-crank PDE (see equation (\ref{eqn:RankCrank}))
 gives a recurrence for the rank moment generating functions that involves
  the triple products
\begin{align}
&C_{2\alpha}(q)C_{2\beta}(q)C_{2\gamma}(q)P(q)^{-2} =: \sum_{n} M_{2\alpha,2\beta,2\gamma}(n) q^{n}
 \end{align}
 for $\alpha, \beta, \gamma\ge 0$.   In this
section we will deal with the asymptotic evaluation of $M_{2\alpha,2\beta,2\gamma}(n).$

\begin{theorem}\label{thm:crankProd}
If $\alpha+\beta+\gamma = k$, then
\begin{align*}
M_{2\alpha,2\beta,2\gamma}(n) =\ & \pi\xi_{2\alpha, 2\beta, 2\gamma}(24 n-1)^{k - 3/4} I_{3/2-2k}(y_{n})+ \xi'_{2\alpha,2\beta,2\gamma} (24n-1)^{k-5/4}I_{3/2-2k+1}(y_{n})\\ & +O\(n^{k-7/4}\cdot n^{-1/4} \exp(y_n) \),
\end{align*}
with
$\xi_{2\alpha, 2\beta, 2\gamma} := \frac{1}{4} \xi_{2\alpha}\xi_{2\beta}\xi_{2\gamma}$
 and
 $$
\xi_{2\alpha, 2\beta, 2\gamma}' :=
\frac14 \left(  \xi_{2\alpha}'\xi_{2\beta}\xi_{2\gamma}
+ \xi_{2\alpha}\xi_{2\beta}'\xi_{2\gamma}
+ \xi_{2\alpha}\xi_{2\beta}\xi_{2\gamma}' \right).
 $$
\end{theorem}
\noindent We again use Lemma \ref{lem:besselFunc} to shift the indices for easier comparisons.
\begin{corollary}
The asymptotic expansion of the triple crank product is given by
\begin{align*}
M_{2\alpha,2\beta,2\gamma}(n) =\ &
\pi\xi_{2\alpha, 2\beta, 2\gamma}
(24 n-1)^{k - 3/4} I_{3/2}(y_{n})+ \widetilde{\xi}_{2\alpha,2\beta,2\gamma} (24n-1)^{k-5/4}I_{1/2}(y_{n})\\ & +O\(n^{k-7/4}\cdot n^{-1/4} \exp(y_n) \),
\end{align*}
where
$$
 \widetilde{\xi}_{2\alpha, 2\beta, 2\gamma}
 :=
-3 (2k) \left(2k-3 \right)  \xi_{2\alpha, 2\beta, 2\gamma}+
\xi_{2\alpha, 2\beta, 2\gamma}' .
$$
\end{corollary}

\begin{proof}[Proof of Theorem \ref{thm:crankProd}]
Using \eqref{eqn:crankREC} for each term in the product of $C_{2\alpha}
C_{2\beta}C_{2\gamma}P^{-2}$ as well as a minor modification of Lemma \ref{lem:Eisen}
we obtain the recursion
\begin{align}\label{eqn:xiProdREC}
\xi_{2\alpha, 2\beta, 2\gamma} = &-2^{3}
\sum_{\begin{subarray}{c}
1\le j \le \alpha \\ 1 \le i \le \beta \\ 1 \le \ell \le \gamma
\end{subarray}} \binom{2\alpha - 1}{2j-1}\binom{2\beta-1}{2i-1}\binom{2\gamma -1}{2\ell-1}
\frac{B_{2j}}{4j} \frac{B_{2i}}{4i} \frac{B_{2\ell}}{4\ell} (-1)^{i+j+\ell} \xi_{2\alpha - 2j, 2\beta - 2i, 2\gamma - 2\ell}.
\end{align}
Induction and three applications of the recurrence (\ref{eqn:xiREC})
 shows that $\xi_{2\alpha, 2\beta, 2\gamma} = \frac{1}{4} \xi_{2\alpha}\xi_{2\beta}\xi_{2\gamma}$ is the unique solution to (\ref{eqn:xiProdREC}).

 We next turn to $\xi_{2\alpha,2\beta,2\gamma}'$.  Following similar arguments as those that led to the
recursion for $\xi_{2k}'$, we obtain
\begin{align*}
\xi_{2\alpha, 2\beta, 2\gamma}'
 =
& 8
\sum_{\begin{subarray}{c}
1\le j \le \alpha \\ 1 \le i \le \beta \\ 1 \le \ell \le \gamma
\end{subarray}} \binom{2\alpha - 1}{2j-1}\binom{2\beta-1}{2i-1}\binom{2\gamma -1}{2\ell-1}
\frac{B_{2j}}{4j} \frac{B_{2i}}{4i} \frac{B_{2\ell}}{4\ell} (-1)^{i+j+\ell+3} \xi_{2\alpha - 2j, 2\beta - 2i, 2\gamma - 2\ell}' \\
& +8  \sum_{\begin{subarray}{c}
 1 \le i \le \beta \\ 1 \le \ell \le \gamma
\end{subarray}} (2\alpha - 1) \binom{2\beta-1}{2i-1}\binom{2\gamma -1}{2\ell-1}
6 \,
\frac{B_{2}}{4} \frac{B_{2i}}{4i} \frac{B_{2\ell}}{4\ell} (-1)^{i+\ell+3} \xi_{2\alpha - 2, 2\beta - 2i, 2\gamma - 2\ell}\\
&+8\sum_{\begin{subarray}{c}
1\le j \le \alpha  \\ 1 \le \ell \le \gamma
\end{subarray}} \binom{2\alpha - 1}{2j-1} (2\beta-1)\binom{2\gamma -1}{2\ell-1}
6 \,
\frac{B_{2j}}{4j} \frac{B_{2}}{4} \frac{B_{2\ell}}{4\ell} (-1)^{j+\ell+3} \xi_{2\alpha - 2j, 2\beta - 2, 2\gamma - 2\ell}\\
&+8 \sum_{\begin{subarray}{c}
1\le j \le \alpha \\ 1 \le i \le \beta
\end{subarray}} \binom{2\alpha - 1}{2j-1}\binom{2\beta-1}{2i-1} (2\gamma -1)
6 \, \frac{B_{2j}}{4j} \frac{B_{2i}}{4i} \frac{B_{2}}{4} (-1)^{i+j+3} \xi_{2\alpha - 2j, 2\beta - 2i, 2\gamma - 2}.
\end{align*}
Again using the formula $\xi_{2\alpha, 2\beta, 2\gamma} = \frac{1}{4} \xi_{2\alpha} \xi_{2\beta} \xi_{2\gamma}$ and  (\ref{eqn:xiREC}), we may simplify the second line (and analogously the third and fourth) as
\begin{align*}
& -\frac{2 (2 \alpha-1)}{16} \xi_{2 \alpha-2}
\left(
2  \sum_{1 \le i \le \beta}  \binom{2\beta-1}{2i-1}
\frac{B_{2i}}{4i}  (-1)^{i+1} \xi_{ 2\beta - 2i}
\right)
\left(
2  \sum_{1 \le \ell \le \gamma}  \binom{2\gamma-1}{2\ell-1}
\frac{B_{2\ell}}{4\ell}  (-1)^{\ell+1} \xi_{ 2\gamma - 2\ell}
\right)\\
 & \qquad \qquad \qquad \qquad = -\frac{2(2\alpha - 1)}{16} \xi_{2\alpha - 2} \xi_{2 \beta} \xi_{2\gamma} = \frac{1}{4\alpha}\xi_{2 \alpha}' \xi_{2\beta}\xi_{2 \gamma}.
\end{align*}
This yields the following recursion for $\xi_{2\alpha,2\beta,2\gamma}'$
\begin{align}\label{eqn:xi'Triple}
\xi_{2\alpha, 2\beta, 2\gamma}' &=
8
\sum_{\begin{subarray}{c}
1\le j \le \alpha \\ 1 \le i \le \beta \\ 1 \le \ell \le \gamma
\end{subarray}} \binom{2\alpha - 1}{2j-1}\binom{2\beta-1}{2i-1}\binom{2\gamma -1}{2\ell-1}
\frac{B_{2j}}{4j} \frac{B_{2i}}{4i} \frac{B_{2\ell}}{4\ell} (-1)^{i+j+\ell+3} \xi_{2\alpha - 2j, 2\beta - 2i, 2\gamma - 2\ell}'   \\
&
+\frac{1}{4\alpha} \xi_{2 \alpha}'\xi_{2 \beta} \xi_{2\gamma}
 +\frac{1}{4\beta} \xi_{2 \alpha}\xi_{2 \beta}' \xi_{2\gamma}
 +\frac{1}{4\gamma} \xi_{2 \alpha}\xi_{2 \beta} \xi_{2\gamma}'. \notag
\end{align}
Our claimed formula $\xi_{2\alpha, 2\beta, 2\gamma}' =
\frac14 (\xi_{2\alpha}'\xi_{2\beta}\xi_{2\gamma}
+ \xi_{2\alpha}\xi_{2\beta}'\xi_{2\gamma}
+ \xi_{2\alpha}\xi_{2\beta}\xi_{2\gamma}')$ solves this recurrence, as using (\ref{eqn:xiREC}) and (\ref{eqn:xi'ID}) to evaluate the right side of (\ref{eqn:xi'Triple}) gives the expected result:
\begin{align*}
\frac{1}{4}\left( 1-\frac{1}{\alpha}\right) \xi_{2 \alpha}'\xi_{2 \beta} \xi_{2\gamma}
& + \frac{1}{4}\left( 1-\frac{1}{\beta}\right) \xi_{2 \alpha}\xi_{2 \beta}' \xi_{2\gamma}
+ \frac{1}{4}\left( 1-\frac{1}{\gamma}\right) \xi_{2 \alpha}\xi_{2 \beta} \xi_{2\gamma}' \\
& + \frac{1}{4\alpha}\xi_{2 \alpha}' \xi_{2\beta}\xi_{2 \gamma}
+ \frac{1}{4\beta}\xi_{2 \alpha} \xi_{2\beta}'\xi_{2 \gamma}
+  \frac{1}{4\gamma}\xi_{2 \alpha}\xi_{2\beta}\xi_{2 \gamma}'
\\
 & =
\frac14 \left(  \xi_{2\alpha}'\xi_{2\beta}\xi_{2\gamma}
+ \xi_{2\alpha}\xi_{2\beta}'\xi_{2\gamma}
+ \xi_{2\alpha}\xi_{2\beta}\xi_{2\gamma}' \right)
=
 \xi_{2\alpha, 2\beta, 2\gamma}'.
\end{align*}

\end{proof}

\section{Rank Asymptotics}\label{sec:Rank}
In this section we prove the asymptotic expansion for $N_{2k}
(n)$ given in Theorem \ref{thm:rank}.  In order to relate the rank moments to the crank moments that we have already calculated, we use Atkin and Garvan's rank-crank PDE, which states that
 \begin{align} \label{eqn:RankCrank}
\sum_{i=0}^{k-1}\binom{2k}{2i} \sum_{\begin{subarray}{c} 2\alpha +2\beta +2 \gamma = 2k-2i\\
\alpha, \beta, \gamma \ge 0  \end{subarray}}&\binom{2k-2i}{2\alpha, 2\beta, 2\gamma} C_{2\alpha}(q)
C_{2\beta}(q)C_{2\gamma}(q)P(q)^{-2} - 3 (2^{2k-1}  - 1) C_{2}(q) \nonumber
\\ =& \frac{1}{2}(2k-1)(2k-2)R_{2k}(q) + 6 \sum_{i=1}^{k-1} \binom{2k}{2i} (2^{2i-1}-1) \delta_{q}(R_{2k-2i}(q))
\\ &+\sum_{i=1}^{k-1} \(\binom{2k}{2i+2} (2^{2i+1}-1) - 2^{2i} \binom{2k}{2i+1} + \binom{2k}{2i} \) R_{2k-2i}(q), \nonumber
\end{align}
where $\delta_{q}: = q \frac{d}{dq}$.
We once more argue inductively in order to find recurrences for $\lambda_{2k}$ and $\widetilde{\lambda}_{2k}$ as defined in Theorem \ref{thm:rank}.

\begin{proof}[Proof of Theorem \ref{thm:rank}]
The rank-crank PDE and Theorems \ref{thm:crank} and \ref{thm:crankProd} inductively imply that there is an asymptotic expansion for the rank moments of the form
$$
N_{2k}(n) = \pi \lambda_{2k} (24n-1)^{k-3/4}I_{3/2-2k}(y_{n}) + \lambda_{2k}'
(24n-1)^{k-5/4}I_{3/2-2k-1}(y_{n}) + O\left(n^{k - 2} e^{y_{n}}\right).
$$
Lemma \ref{lem:besselFunc} then implies that
$$
N_{2k}(n) = \pi \lambda_{2k} (24n-1)^{3/4}I_{3/2}(y_{n}) +
\widetilde{\lambda}_{2k}
(24n-1)^{k-5/4}I_{1/2}(y_{n}) + O\left(n^{k - 2} e^{y_{n}}\right),
$$
with $\widetilde{\lambda}_{2k}:= -3 (2k) (2k-3)\lambda_{2k}  + \lambda_{2k}'.$

The only terms that contribute to the leading order on the crank side of (\ref{eqn:RankCrank}) are
\begin{equation}\label{eqn:PDEcrank}
\sum_{\alpha + \beta + \gamma = k} \binom{2k}{2\alpha, 2\beta, 2\gamma} \xi_{2\alpha, 2\beta, 2\gamma}
\end{equation}
 while the only terms contributing to the main term on the rank side of the equation are
\begin{equation}\label{eqn:PDErank}
 \binom{2k-1}{2} \lambda_{2k} + \binom{2k}{2} \frac{1}{4} \lambda_{2k-2}.
\end{equation}
We can now show inductively that $\lambda_{2k-2} = \xi_{2k-2}$.  Indeed, adding the terms from (\ref{eqn:PDEcrank}) and (\ref{eqn:PDErank}) and using the formula for $\xi_{2\alpha, 2\beta, 2\gamma}$ shows that the claimed equality is equivalent to the triple product summation
$$\sum_{\alpha, \beta, \gamma} \binom{2k}{2\alpha, 2\beta, 2\gamma} B_{2\alpha}(1/2)
B_{2\beta}(1/2) B_{2\gamma}(1/2) = \binom{2k-1}{2} B_{2k}(1/2) - \binom{2k}{2} \frac{1}{4} B_{2k-2}
(1/2).$$
Like (\ref{eqn:bern}), this is also a specialization of a general Bernoulli polynomial identity that can be found in \cite{Dil}.  As a consequence of the equality between $\xi_{2k}$ and $\lambda_{2k}$, we now have
\begin{equation}\label{eqn:xiTriple}
\sum_{\alpha + \beta + \gamma = k} \binom{2k}{2\alpha, 2\beta, 2\gamma} \xi_{2\alpha, 2\beta, 2\gamma} =
\binom{2k-1}{2} \xi_{2k} + \binom{2k}{2} \frac{1}{4} \,\xi_{2k-2}.
\end{equation}

We next consider the second leading term.
Using similar reasoning as above, we equate the terms of second highest order from (\ref{eqn:RankCrank}), obtaining
\begin{equation} \label{eqn:secondterm}
\sum_{\alpha + \beta + \gamma = k} \binom{2k}{2\alpha, 2\beta, 2\gamma}
\widetilde{\xi}_{2\alpha, 2\beta, 2\gamma}
=
\binom{2k-1}{2} \widetilde{\lambda}_{2k} + \binom{2k}{2} \frac{1}{4} \widetilde{\lambda}_{2k-2}.
\end{equation}
Expanding the left side of (\ref{eqn:secondterm}) gives
\begin{equation}\label{eqn:left}
\sum_{\alpha + \beta + \gamma = k} \binom{2k}{2\alpha, 2\beta, 2\gamma}
\widetilde{\xi}_{2\alpha, 2\beta, 2\gamma}
=
\sum_{\alpha + \beta + \gamma = k} \binom{2k}{2\alpha, 2\beta, 2\gamma}
\left(
-3\cdot 2k  \left(2k-3\right) \xi_{2\alpha, 2\beta, 2\gamma} +
\xi_{2\alpha, 2\beta, 2\gamma}'
\right).
\end{equation}
The first term of (\ref{eqn:left}) can be evaluated using (\ref{eqn:xiTriple}):
\begin{align}\label{eqn:left1}
-3\cdot 2k  \left(2k-3\right)
& \sum_{\alpha + \beta + \gamma = k} \binom{2k}{2\alpha, 2\beta, 2\gamma}
 \xi_{2\alpha, 2\beta, 2\gamma}
 =
 -3\cdot 2k  \left(2k-3\right)
 \left(
 \binom{2k-1}{2} \xi_{2k} + \binom{2k}{2} \frac{1}{4} \xi_{2k-2}
 \right)  \notag \\
  & = -\frac32 (2k)(2k-1)(2k-2)(2k-3)\xi_{2k}
  - \frac38(2k)^2(2k-1)(2k-3)\xi_{2k-2},
 \end{align}
and the second term of (\ref{eqn:left}) is
 \begin{align}\label{eqn:left2}
 \sum_{\alpha + \beta + \gamma = k} & \binom{2k}{2\alpha, 2\beta, 2\gamma}
\xi_{2\alpha, 2\beta, 2\gamma}'
= \frac14  \sum_{\alpha + \beta + \gamma = k} \binom{2k}{2\alpha, 2\beta, 2\gamma}
 \left(  \xi_{2\alpha}'\xi_{2\beta}\xi_{2\gamma}
+ \xi_{2\alpha}\xi_{2\beta}'\xi_{2\gamma}
+ \xi_{2\alpha}\xi_{2\beta}\xi_{2\gamma}' \right) \notag
   \\
& =
 \frac14  \sum_{\alpha + \beta + \gamma = k} \binom{2k}{2\alpha, 2\beta, 2\gamma}
 \left(
 -\frac14 (2 \alpha)(2 \alpha-1)
  \xi_{2\alpha-2}\xi_{2\beta}\xi_{2\gamma}
  -\frac14 (2 \beta)(2 \beta-1)
  \xi_{2\alpha}\xi_{2\beta-2}\xi_{2\gamma} \right. \notag \\
  & \left. \qquad \qquad \qquad \qquad \qquad \qquad \qquad \qquad \qquad \qquad \qquad \qquad
   -\frac14 (2 \gamma)(2 \gamma-1)
    \xi_{2\alpha}\xi_{2\beta}\xi_{2\gamma-2} \right).
\end{align}
After using the identity
$\binom{2k}{2\alpha, 2\beta, 2\gamma}(2\alpha)(2\alpha - 1) = 2k(2k-1)\binom{2k-2}{2\alpha-2, 2\beta, 2\gamma}$ (with similar analogues for $\beta$ and $\gamma$) and making a shift of indices, each of the three terms of (\ref{eqn:left2}) become alike. Using the definition of  $\xi_{2\alpha,2\beta,2\gamma}$ this gives  a total of
\begin{align}\label{eqn:left2Final}
-\frac34 & (2k)(2k-1)
   \sum_{\alpha + \beta + \gamma = k-1} \binom{2k-2}{2\alpha, 2\beta, 2\gamma}
  \xi_{2\alpha,2\beta,2\gamma}
    \\
    & =
    -\frac{3}{4} (2k)(2k-1)
     \left(
 \binom{2k-3}{2} \xi_{2k-2} + \binom{2k-2}{2} \frac{1}{4} \xi_{2k-4}
 \right) \notag
 \\
 & = -\frac{3}{8} (2k)(2k-1)(2k-3)(2k-4)\xi_{2k-2}
 -\frac{3}{32}(2k)(2k-1)(2k-2)(2k-3)\xi_{2k-4}. \notag
 \end{align}
 Combining (\ref{eqn:left1}) and (\ref{eqn:left2Final}) gives
 \begin{multline*}
 \sum_{\alpha + \beta + \gamma = k} \binom{2k}{2\alpha, 2\beta, 2\gamma}
\widetilde{\xi}_{2\alpha, 2\beta, 2\gamma}
=  -\frac32 (2k)(2k-1)(2k-2)(2k-3)\xi_{2k} \\
-\frac{3}{4} (2k)(2k-1)(2k-2)(2k-3)  \xi_{2k-2}
 -\frac{3}{32}(2k)(2k-1)(2k-2)(2k-3)\xi_{2k-4}.
 \end{multline*}

It is now easy to see that the expression
 $$
 \widetilde{\lambda}_{2k}
= -3\cdot 2k(2k-3) \xi_{2k}
 -\frac34 2k(2k-1) \xi_{2k-2}
 $$
from the theorem statement is the unique solution to (\ref{eqn:secondterm}).
\end{proof}

\end{document}